\documentclass[reqno]{amsart}

\usepackage{amsmath,amssymb,amsthm,subcaption}

\usepackage{tikz}
\usetikzlibrary{calc,intersections,through,backgrounds}

\usepackage{caption}
\usepackage{graphicx}
\usepackage{graphics}

\newtheorem{theorem}{Theorem}
\newtheorem*{thm}{Theorem}

\newtheorem*{corollary}{Corollary}

\begin{document}

\title[]{Magnetic Schr\"odinger Operators \\and Landscape Functions}

\author[]{Jeremy G. Hoskins}
\address{Department of Statistics, University of Chicago,
IL 60637} \email{ jeremyhoskins@uchicago.edu}

\author[]{Hadrian Quan}
\address{Department of Mathematics, University of Washington, Seattle, WA 98195, USA} \email{hadrianq@uw.edu}

\author[]{Stefan Steinerberger}
\address{Department of Mathematics, University of Washington, Seattle, WA 98195, USA} \email{steinerb@uw.edu}

\thanks{S.S. is supported by the NSF (DMS-2123224) and the Alfred P. Sloan Foundation.}

\begin{abstract} We study localization properties of low-lying eigenfunctions of magnetic Schr\"odinger operators 
 $$\frac{1}{2} \left(- i\nabla - A(x)\right)^2 \phi + V(x) \phi = \lambda \phi,$$
where $V:\Omega \rightarrow \mathbb{R}_{\geq 0}$ is a given potential and $A:\Omega \rightarrow \mathbb{R}^d$ induces a magnetic field.
We extend the Filoche-Mayboroda inequality and prove a refined inequality in the magnetic setting which can predict the points where low-energy eigenfunctions are localized. This result is new even in the case of vanishing magnetic field. Numerical examples illustrate the results.\end{abstract}
\maketitle

\section{Introduction}
\subsection{Introduction.}
A fundamental problem in mathematics and physics is to understand the behavior of low-energy eigenstates of Schr\"odinger operators. Given an open and bounded domain $\Omega \subset \mathbb{R}^d,$ and a potential $V: \Omega \rightarrow \mathbb{R}_{\geq 0}$, this amounts to characterizing solutions of the equation
\begin{align*}
-\Delta \phi(x) + V(x) \phi(x) &= \lambda \phi(x)~ \quad \mbox{in~}\Omega \\
 \phi&= 0 \qquad \quad \mbox{on}~  \partial \Omega
\end{align*}
for small values of $\lambda$. Our primary focus in this paper is on the regime where $V$ oscillates rapidly at small scales. It is known that solutions of this equation may become strongly localized \cite{anderson}. Moreover, it is also understood that in this setting the boundary conditions are not tremendously important. We mention that in the case $V=0$ it is possible to have fascinating localization phenomena that stem from Neumann boundary conditions (see \cite{sapo1, jones,  sapo2}) but these will not be discussed here.  Our main question is whether it is possible to predict efficiently, using only $V$,
where such highly-localized eigenfunctions will be concentrated. 

\subsection{Landscape.} This question has received renewed attention in recent years after the introduction of the landscape function by Filoche \& Mayboroda \cite{fil}.
They define the landscape function as $u:\Omega: \mathbb{R} \rightarrow \mathbb{R}_{}$ satisfying
\begin{align*}
(-\Delta + V)u &=1~ \qquad \mbox{in~}\Omega \\
 u&= 0 \qquad  \mbox{on}~~ \partial \Omega
\end{align*}
and prove that $u$ exerts pointwise control on all eigenfunctions $(-\Delta + V)\phi = \lambda \phi$ 
$$ |\phi(x)|  \leq \lambda u(x)  \|\phi\|_{L^{\infty}(\Omega)} .$$ 
The landscape function turns out to have remarkable predictive power: the first few eigenfunctions tend to localize close to the maxima of $u$. Indeed, one can think of $1/u$ as a suitable regularization of the potential $V$ with the property that eigenfunctions localize in its local minima. This has inspired considerable subsequent work, we refer to \cite{arnold0, arnold, arnold2, alt, altmann, chal, david, fil, fil2, fil3, har, korn, lef, leite, lu, pic}.

\subsection{Local Landscape.}  Another approach was proposed by the third author in \cite{steini, steini3}: there it was shown that, on sufficiently small scales, there exists a canonical smoothing of the potential $V$ given by the Wiener integral
$$ V_t(x) =   \mathbb{E} ~\frac{1}{t} \int_0^t V(\omega_x(s)) ds,$$
where $\omega_x(s)$ is a Brownian motion started at $x$ at time $s$. As is shown in \cite{steini3}, this mollifier arises naturally when looking for local approximations of the identity with good error estimates. Moreover, there is a probability-free way of computing this object.
Since Brownian motion is, in free space, distributed according to a Gaussian and since expectation is linear, we can express $V_t \sim V*k_t$, where
$$ k_t(x) = \frac{1}{t} \int_0^t \frac{ \exp\left( - \|x\|^2/ (4s) \right)}{(4 \pi s )^{d/2}} ds$$
with an exponentially small error provided that $t \ll \mbox{dist}(x, \Omega^c)^2$. Note that convolution with $k_t$ can be performed rapidly using the Fast Fourier Transform. The local landscape is, in a suitable sense, the canonical mollifier for small $t$ (see \cite{steini3}) and was empirically shown to have
predictive power roughly comparable to that of the landscape function \cite{li2} for larger values of $t$.

\section{Magnetic Schr\"odinger Operators}
Our goal is to extend and unify these two approaches for magnetic Schr\"odinger operators. 
 Let $\Omega \subset \mathbb{R}^d$ be a bounded domain with smooth boundary (mainly for ease of exposition, this could be relaxed). We consider the eigenvalue problem
$$ \frac{1}{2} \left(- i\nabla - A(x)\right)^2 \phi(x) + V(x)\phi(x) = \lambda \phi(x)$$
where $V:\Omega \rightarrow \mathbb{R}_{\geq 0}$,  $A:\Omega \rightarrow \mathbb{R}^d$ and $\phi$ is subject to Dirichlet boundary conditions. 
The main technical ingredient in our approach is a Feynman-Kac formula for magnetic Schr\"odinger operators due to  Broderix, Hundertmark \& Leschke \cite{broderix}. In particular, we will work with the regularity assumptions $A \in C^1(\mathbb{R}^d, \mathbb{R}^d)$ and $V \in C^1(\mathbb{R}^d, \mathbb{R}_{\geq 0})$ which could be somewhat weakened (see \cite{broderix} for a discussion). In our setting of interest where $V$ is wildly oscillatory, one does not expect the boundary
conditions to have a significant impact on the behavior of eigenfunctions inside the domain.
 The only prior work in this direction 
is due to Poggi \cite{poggi} who proposes a different extension of the landscape function leading to different results.

\subsection{A Magnetic Filoche-Mayboroda inequality} We start by extending the Filoche-Mayboroda Landscape function to the magnetic setting.
We consider the solution of the equation,  henceforth referred to as the landscape function,
 \begin{align*}
 \left[- \frac12 \Delta + V\right]u &= 1 \quad \mbox{in}~\Omega \\
 u &= 0 \quad \mbox{on}~\partial \Omega.
 \end{align*}
 The factor $1/2$ in front of the Laplacian is equivalent to the usual Laplacian up to scaling $V$ and $\lambda$ by a factor of 2: it is the probabilistic normalization of the Laplacian chosen to facilitate comparison with the result in \cite{broderix}.
\begin{thm}[Filoche-Mayboroda \cite{fil}]
For any eigenfunction $(-(1/2) \Delta + V) \phi = \lambda \phi$ subject to Dirichlet boundary conditions and all $x \in \Omega$
$$ \frac{|\phi(x)|}{\| \phi\|_{L^{\infty}}} \leq \lambda \cdot u(x).$$
\end{thm}
Note that the quantity on the left-hand side always assumes the value 1 at the global extrema of any eigenfunction. Therefore,
eigenfunctions corresponding to eigenvalue $\lambda$ can only localize at points $x \in \Omega$ where $u(x) \cdot \lambda \geq 1$. Our first main result shows that this inequality extends to magnetic Schr\"odinger operators. 
 
 \begin{theorem} Let $\Omega \subset \mathbb{R}^d$ be a bounded domain with smooth boundary, let $A \in C^1(\Omega, \mathbb{R}^d)$, let $V \in C^1(\Omega, \mathbb{R}_{\geq 0})$ and let $\phi$ be a solution of
\begin{equation} \label{eq:equa}
 \frac{1}{2} \left(- i\nabla - A(x)\right)^2 \phi(x) + V(x)\phi(x) = \lambda \phi(x)
 \end{equation}
 subject to Dirichlet boundary conditions. Then, for all $x \in \Omega$,
$$ \frac{|\phi(x)|}{\| \phi\|_{L^{\infty}}} \leq \lambda \cdot u(x).$$
\end{theorem}

 An interesting aspect of the inequality is that the magnetic potential $A$ does not explicitly appear in the inequality (it appears implicitly in $\lambda$). Indeed, if $A$ is much smoother than the potential $V$ (in the sense of being virtually constant on the scale of the eigenfunction), this is not particularly surprising since, intuitively, it contributes only a phase modulation to the corresponding eigenfunction (see \S 3).

\subsection{A Refined Localization Inequality} 
We also prove a refined localization inequality which is stronger than Theorem 1 (and will imply Theorem 1).
\begin{theorem} Under the same assumptions as Theorem 1, for all $x \in \Omega$ and $t > 0$,
$$  \frac{|\psi(x_0)|}{\|\psi\|_{L^{\infty}}} \leq   \mathbb{E}_{\omega}\left(\chi_{\Omega}(\omega_{x_0}, t) \exp\left(  \lambda \cdot t  - \int_0^t V(\omega_{x_0}(s)) ds \right)  \right)$$
where expectation is taken with respect to all Brownian motions $\omega_x(s)$ started at $\omega_x(0) = x$ running for $t$ units of time and
$$\chi_{\Omega}(\omega,t)  = \begin{cases} 1 \qquad &\mbox{if}~\forall ~0\leq s \leq t: \omega(s) \in \Omega \\ 0 \qquad &\mbox{otherwise.} \end{cases}$$
\end{theorem}

 Note that when $t$ is sufficiently small compared to the distance $d(x_0, \partial \Omega)$ between $x_0$ and the boundary of the domain (in particular, when $t \ll d(x_0, \partial \Omega)^2$), then virtually all Brownian motions will stay inside the domain for $t$ units of time: the difference between having the term $\chi_{\Omega}(\omega,t)$ and omitting it is exponentially small. One important fact is that Theorem 2 contains all the information contained in Theorem 1 and more: this is made precise in the following Corollary.

\begin{corollary} Under the same assumption as Theorem 1, we have
$$ \inf_{t > 0} ~\mathbb{E}_{\omega}\left( \chi_{\Omega}(\omega_{x_0}, t) \exp\left(\lambda \cdot t - \int_0^t V(\omega_{x_0}(s)) ds \right)  \right) \leq \lambda \cdot u(x).$$
\end{corollary}

We refer to \S 3 where this is illustrated on concrete examples.

 \subsection{Regularized Potential.} The second main implication of Theorem 2 is that it sheds some light on the local landscape function: these are introduced as convolutions of the potential and derived from the short-time behavior of the associated parabolic equation. Here, the local landscape naturally arises in a different manner as follows: there are natural cases when the path integral is relatively concentrated around its mean. In such cases, one would expect that
\begin{align} \label{eq:almost}
  \mathbb{E}_{\omega}\left( \chi_{\Omega}(\omega_{x_0}, t) \exp\left( - \int_0^t V(\omega_{x_0}(s)) ds \right)  \right)  \sim \exp\left( -  \mathbb{E}_{\omega} \int_0^t V(\omega_{x_0}(s)) ds \right)  
  \end{align}
 and that quantity is determined by its leading order term which is exactly the local landscape function since
 $$  \mathbb{E}_{\omega} \int_0^t V(\omega_{x_0}(s)) ds = t \cdot V_t(x) + \mbox{small errors},$$
 where the errors come from particles exiting the domain within $t$ units of time, these errors are exponentially small inside the bulk of the domain as long as $t \lesssim d(x, \partial \Omega)^2$.
 We quickly note yet another relationship between the landscape function and
 the local landscape which follows as a byproduct. Considering Theorem 2 and \eqref{eq:almost}, in order to obtain reasonable bounds, one would hope to find
 a parameter $t$ such that 
 $$ \mathbb{E}_{\omega}~  \frac{1}{t} \int_0^t V(\omega_{x_0}(s)) ds \qquad \mbox{is large.}$$
Since this integral is an average over $V$, the values it can assume are naturally limited by the behavior of $V$. It is by now well understood that the inverse landscape $1/u$ can be interpreted as a suitable regularization of the potential.  We show that the local landscape leads to results that are at least comparable to $1/u(x)$ up to exponentially small factors depending only on the distance to the boundary.
 
 \begin{theorem} Under the same assumptions as Theorem 1, there exists a constant $c_d > 0$ depending only on the dimension such that
 $$  \sup_{t > 0}  ~\mathbb{E} \frac{1}{t} \int_0^t V(\omega_x(s)) ds \geq \frac{1}{u(x) + c_d \cdot \exp\left( -\sqrt{\inf_x V(x) } \cdot d(x, \partial \Omega) \right)}.$$
 \end{theorem}
 This may be interpreted as yet another quantification that the path integral can capture at least as much information as $1/u(x)$. Given the recurring idea that one is the average of the other, this is now perhaps less surprising.

 %

\section{Examples}
We illustrates the results for some explicit examples: the case of $[0,1]$ is done in closed form in \S 3.1, we provide more relevant numerical examples in \S 3.2.
\subsection{An analytical example} We start with a very simple example for which analytic expressions are readily obtainable for many of the quantities of interest, to wit, we consider the equation $-\Delta u = \lambda u$ on $[0,\pi]$ with Dirichlet boundary conditions. It is straightforward to show that the first eigenfunction is
$$ \phi(x) = \sin{(x)} \qquad \mbox{corresponding to an eigenvalue of} \quad \lambda = 1.$$
The Filoche-Mayboroda landscape function is $u(x) = x(\pi-x)/2$. We now illustrate Theorem 2 on the same example. Observe that, since $V \equiv 0$,
$$  \mathbb{E}_{\omega}\left( \chi_{\Omega}(\omega_{x_0}, t) \exp\left(\lambda \cdot t - \int_0^t V(\omega_{x_0}(s)) ds \right)  \right)= e^t  \cdot \mathbb{E}_{\omega}  \chi_{\Omega}(\omega_{x_0}, t).$$

\begin{center}
\begin{figure}[h!]
\includegraphics[width=0.6\textwidth]{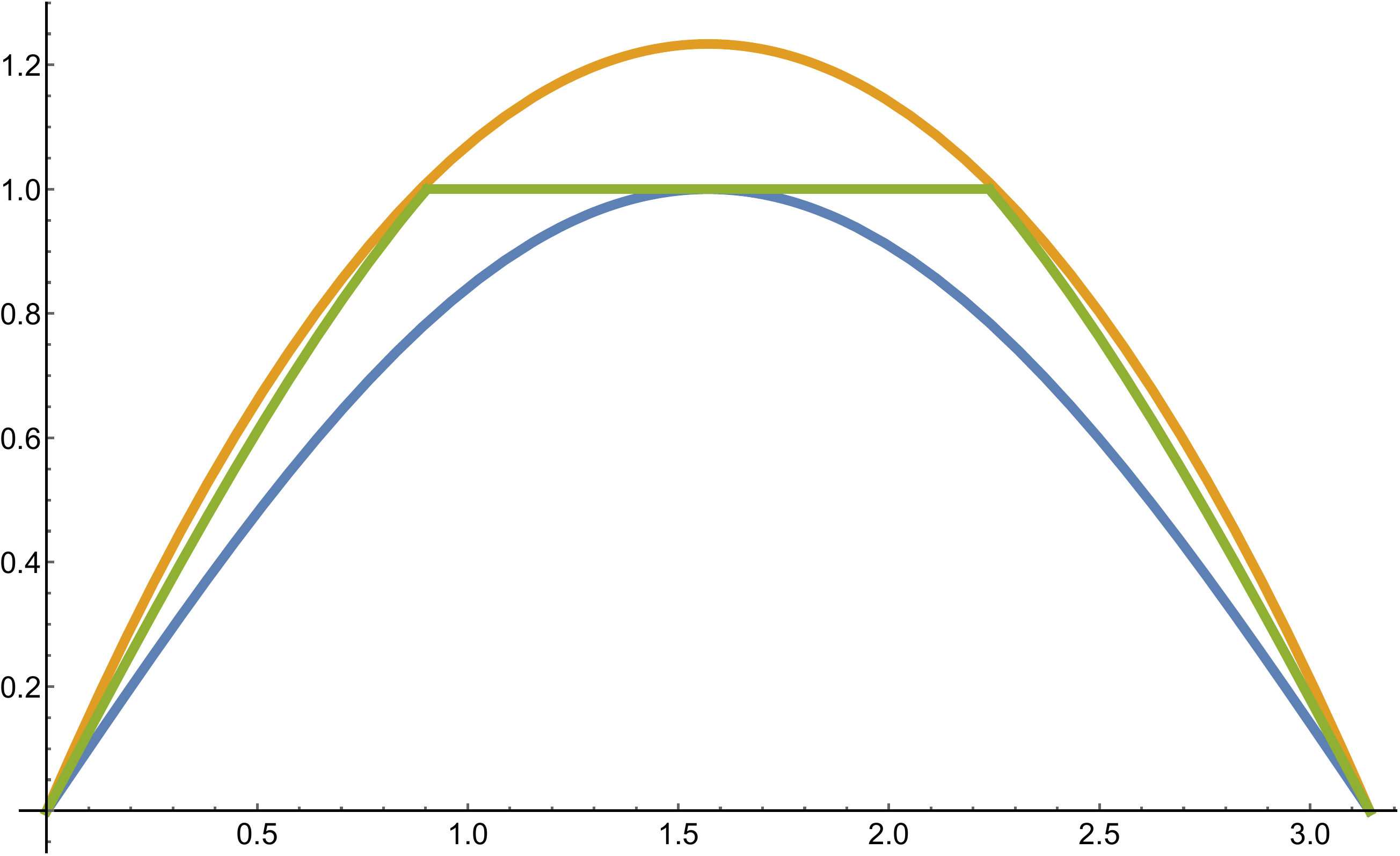}
\caption{The profile of $\sin{x}$ (in blue), the landscape function (in orange) and the profile guaranteed by Theorem 2 (in green).}
\end{figure}
\end{center}
\vspace{-10pt}

This is the survival likelihood of a particle started at $x_0$ and run for $t$ units of time. Next, we note that the distribution of particles started at $x_0$ and run for $t$ units of time while never hitting the boundary is exactly the heat kernel
$$ p_t(x_0, x) =  \frac{2}{\pi} \sum_{k=1}^{\infty} e^{-k^2 t} \sin{(k x_0)} \sin{(k x)}.$$
and therefore
\begin{align*}
e^t  \cdot \mathbb{E}_{\omega}  \chi_{\Omega}(\omega_{x_0}, t) &= e^t \int_0^{\pi} \frac{2}{\pi} \sum_{k=1}^{\infty} e^{-k^2 t} \sin{(k x_0)} \sin{(k x)}\; dx\\
&=  \frac{4 e^t}{\pi} \sum_{k=1 \atop k~\mbox{\tiny odd}}^{\infty} e^{-k^2 t} \frac{\sin{(k x_0)}}{k}. 
\end{align*}
Note that this inequality is true for all $t>0$, allowing one to optimize over $t$. Though for a general $t$, the previous expression is somewhat unwieldy, in the small and large $t$ limits one can obtain relatively compact expressions: 
taking $t \rightarrow 0$
$$ \lim_{t \rightarrow 0} e^t  \cdot \mathbb{E}_{\omega}  \chi_{\Omega}(\omega_{x_0}, t) = 1,$$
since a characteristic function is at most 1.
Simultaneously,
$$ \lim_{t \rightarrow \infty}  \frac{4 e^t}{\pi} \sum_{k=1 \atop k~\mbox{\tiny odd}}^{\infty} e^{-k^2 t} \frac{\sin{(k x_0)}}{k} = \frac{4}{\pi} \sin{(x)}.$$
Therefore, Theorem 2 and the Corollary guarantee that
$$ \sin{(x)} = \frac{|\phi(x)|}{\| \phi\|_{L^{\infty}}} \leq \min \left\{1, \frac{4}{\pi} \sin{(x)} \right\} \leq \frac{x (\pi-x)}{2} = u(x).$$
The three profiles are shown in Figure 1.

\begin{center}
    \begin{figure}[h!]
        \centering
      \includegraphics[width=0.6\textwidth]{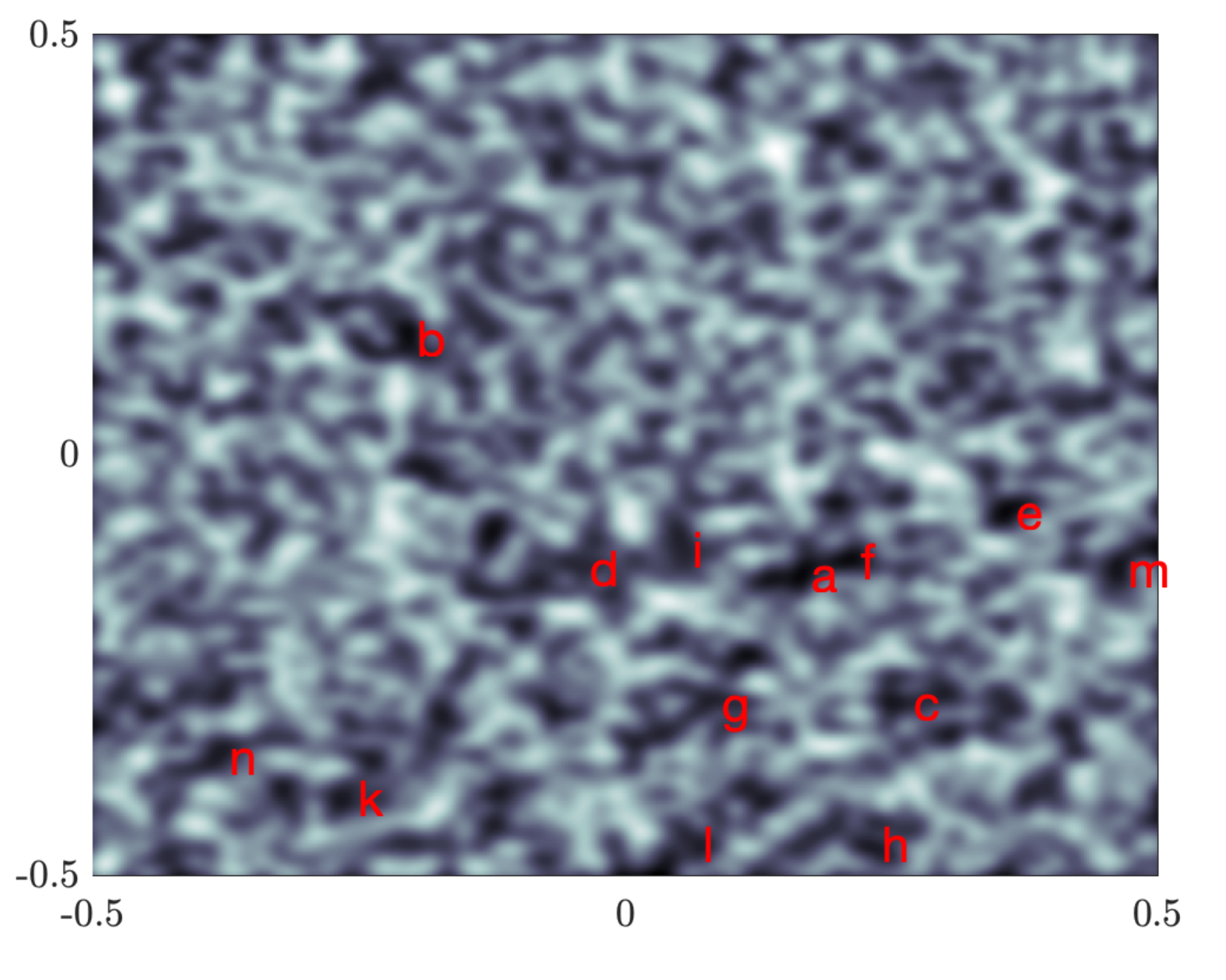}
      \vspace{-10pt}
\caption{The potential $V$ used in the numerical experiments. The approximate locations of the first 13 eigenfunctions (for the case of vanishing magnetic field $\vec{A}\equiv 0$) are given by the red letters.}
    \end{figure}
\end{center}

\vspace{-20pt}

 \subsection{Numerical examples.}
 
This section provides numerical experiments illustrating the results of Theorems 1 and 2 in two dimensions. The domain under consideration is $\Omega = [-1/2,1/2]^2$ with Dirichlet boundary conditions. We constructed an explicit $V$ as a superposition of random Gaussians via
 $$V(x,y) = \sum_{i_x,i_y=0}^{63} v_{i_x,i_y} \exp\left(-\frac{(x-i_x/63+1/2)^2+(x-i_y/63+1/2)^2}{ \sigma^2}\right)$$
 where the weights were i.i.d. uniform random variables on the interval $[0,20480]$ and $\sigma = 1/\sqrt{64}.$ The particular realization that was used is shown in Figure 2. The components of the magnetic potential $A$ were constructed via
 $$\vec{A}_x(x,y) = b\sum_{i_x,i_y=0}^{31} a^{(x)}_{i_x,i_y} \exp\left(-\frac{(x-i_x/31+1/2)^2+(x-i_y/31+1/2)^2}{ \sigma_x^2}\right)$$
  $$\vec{A}_y(x,y) = -b\sum_{i_x,i_y=0}^{31} a^{(y)}_{i_x,i_y} \exp\left(-\frac{(x-i_x/31+1/2)^2+(x-i_y/31+1/2)^2}{ \sigma_y^2}\right)$$
 where $b$ is an arbitrary real number corresponding to the strength of the magnetic field, the weights were i.i.d. normal random variables, and $\sigma_x = \sigma_y = 1/\sqrt{32}.$ In our experiments, we considered $b$ in the range $[0,140].$ The components of the magnetic potential are shown in Figure 3.
 \begin{figure}[h!]
\begin{subfigure}{.5\textwidth}
  \centering
  \includegraphics[width=.99\linewidth]{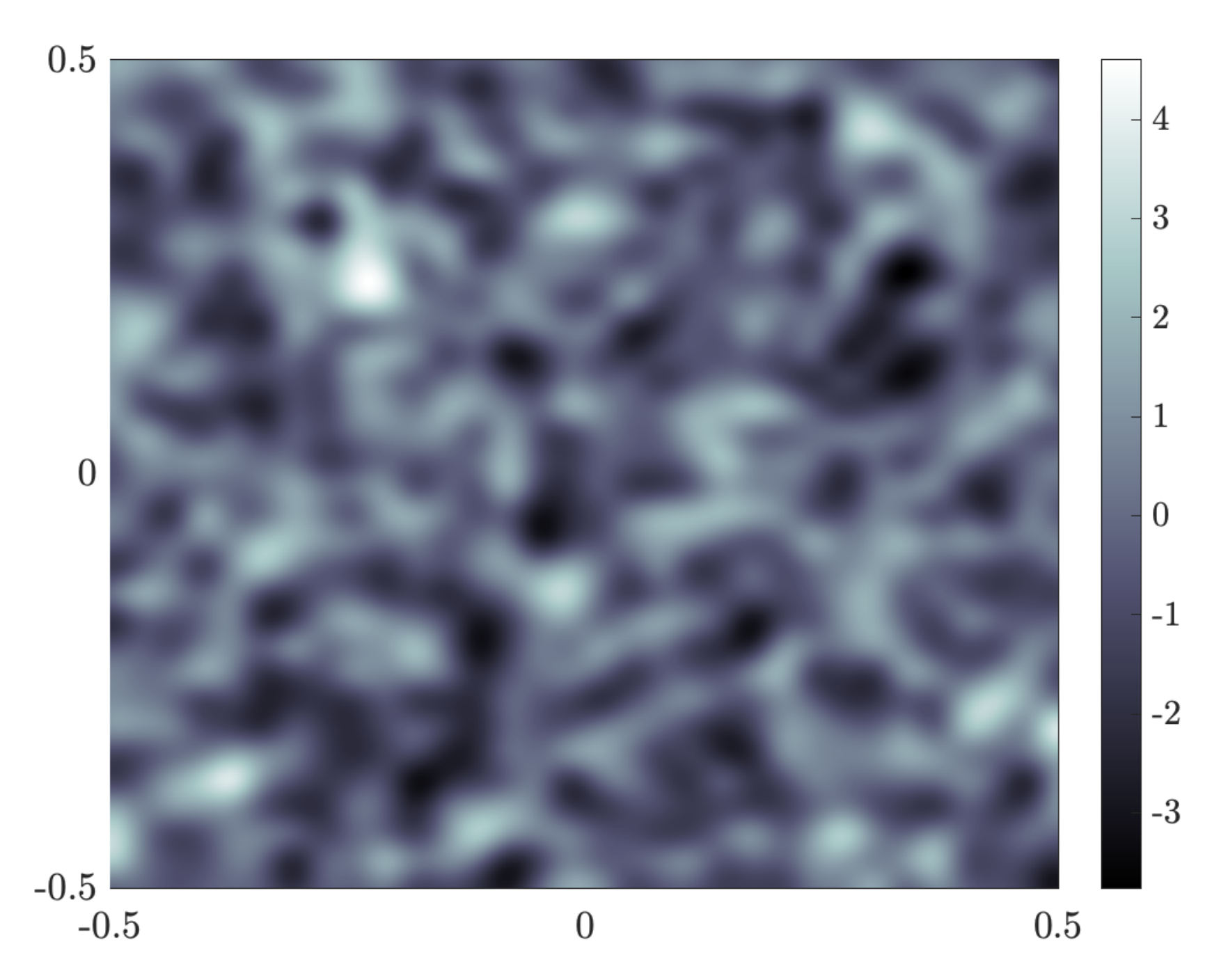}
  \label{fig:magx}
\end{subfigure}%
\begin{subfigure}{.5\textwidth}
  \centering
  \includegraphics[width=.99\linewidth]{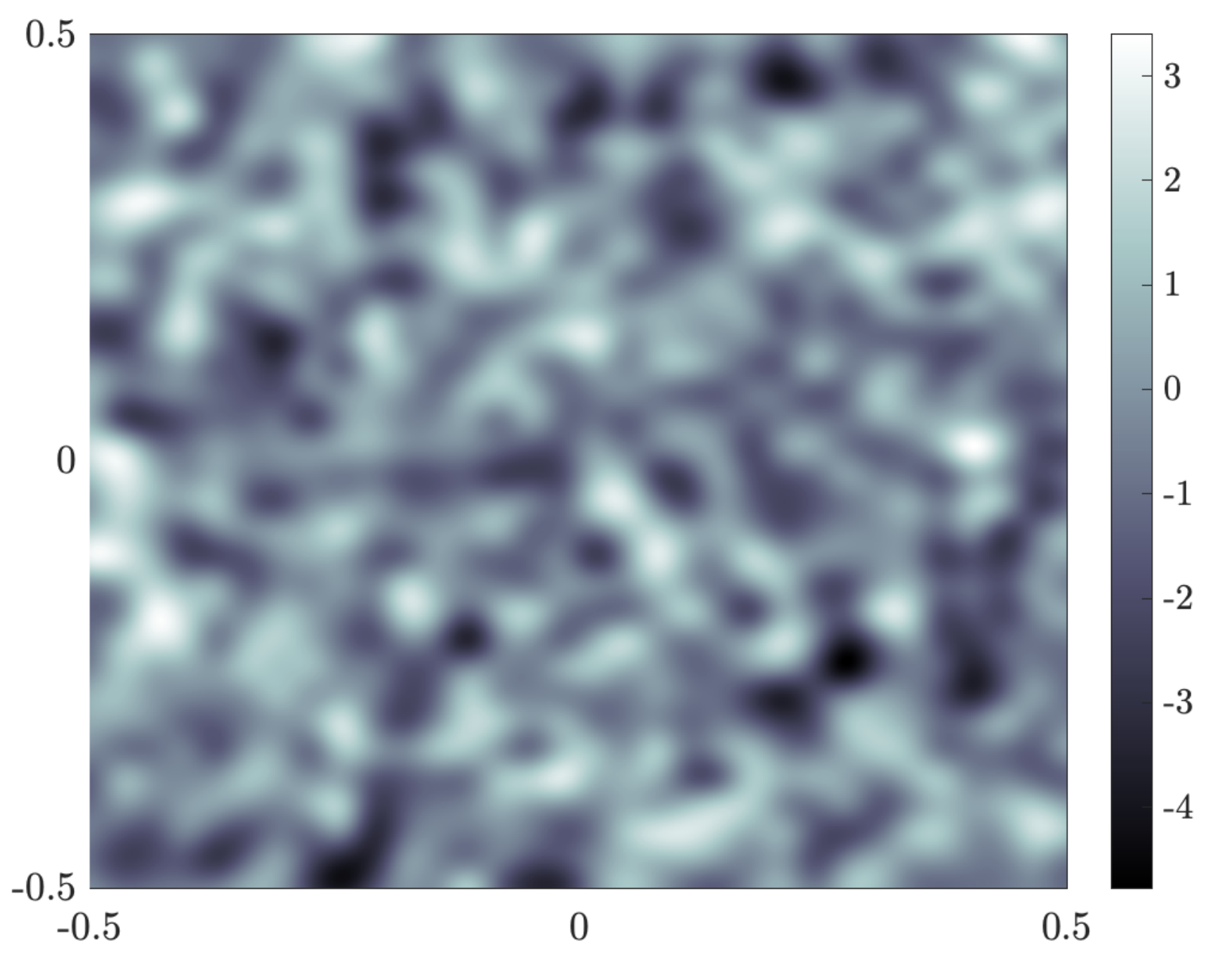}
\label{fig:magy}
  \end{subfigure}
        \vspace{-30pt}
  \caption{The $x-$coordinate (left) and $y-$coordinate (right) of the magnetic potential $A$ that was used for the example.}
 \end{figure}
The landscape function was obtained using a hierarchical Poincar\'{e}-Steklov (HPS) method. The algorithm employed here, based on minor modifications of \cite{gillman}, uses a multidomain spectral collocation discretization method to discretize the problem and a nested-dissection type direct solver to compute a compressed representation of the inverse. We refer the interested reader to the papers \cite{fortunato,gillman2014,gillman} as well as the book \cite{martinsson_book} for details on the method, related algorithms, and other applications. In our experiments, we subdivided the square into $128\times 128$ equally-sized squares and used a $16$th order spectral method on each square for the local solves.
  In order to find the eigenfunctions and eigenvalues, we first used an HPS method to evaluate the `monitor' function $\phi_{f,g}$ defined by
 $$\phi_{f,g}(\lambda) = \left\{\int_\Omega f(x,y) \left[\left(-i\nabla -A\right)^2+V - \lambda\right]^{-1} (g)(x,y)\,{\rm d}x\,{\rm d}y\right\}^{-1}$$ 
 for a range of $\lambda.$ Here $(\left(-i\nabla -A\right)^2+V - \lambda)^{-1}$ denotes the solution operator for Dirichlet boundary conditions, and $f$ and $g$ are arbitrary smooth complex-valued functions. Local minima of $|\phi|$ were used as initializations for inverse power method (the iterates of which were also computed using the HPS solver).\\
 
  \vspace{-10pt}
 \begin{center}
\begin{figure}[h!]
\includegraphics[width=0.7\textwidth]{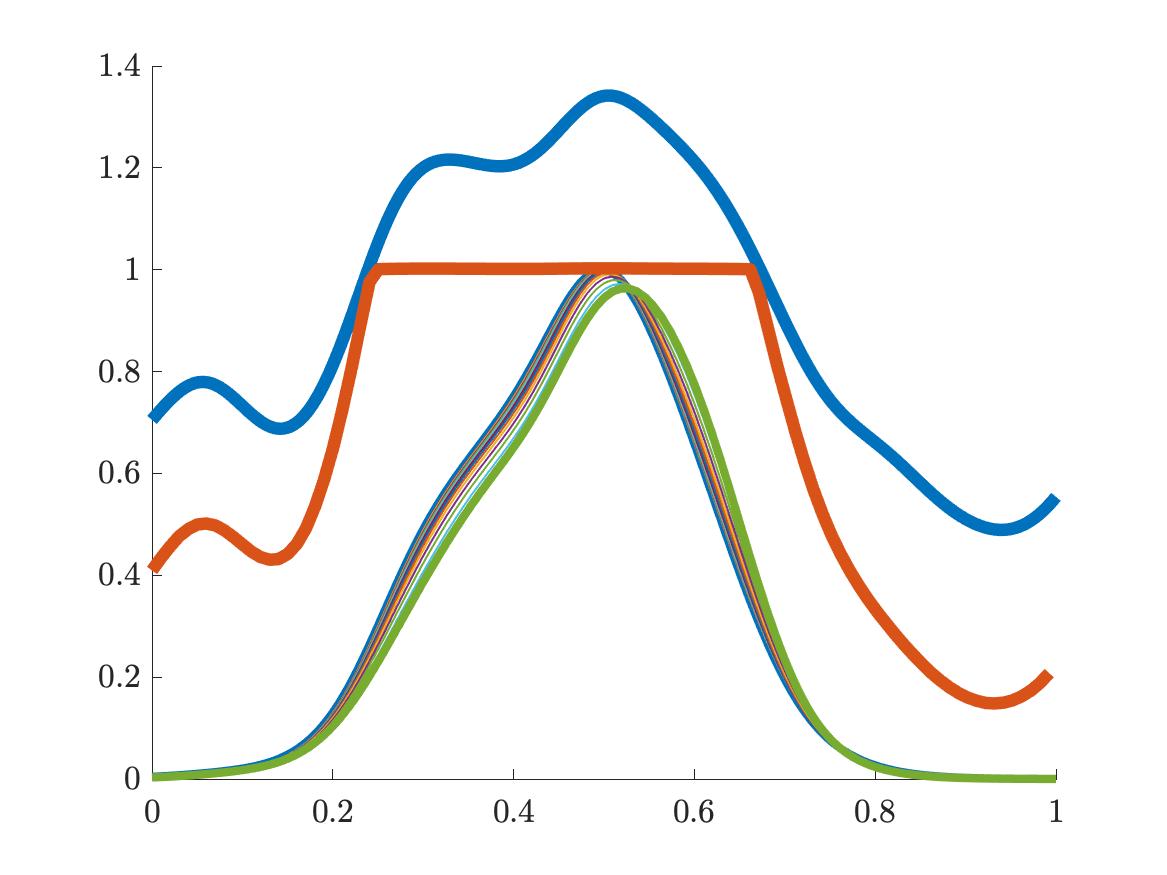}
\caption{The first eigenfunction along a one-dimensional cross-section for different $0 \leq b \leq 140$, the Filoche-Mayboroda landscape function (blue) and the bound coming from Theorem 2 (red).}
\end{figure}
\end{center}
 \vspace{-10pt}
 
We consider the first eigenfunction for the potential and magnetic field shown in Fig. 2 and Fig. 3.
 Fig. 4 shows the first eigenfunction along a one-dimensional line segment centered around its point of localization (a line segment from $(0.05,-0.19)$ to $(0.27,-0.08)$). Fig. 4 shows how the absolute value of the first eigenfunction evolves as a parameter of the strength $b$ of the magnetic potential: as already hinted at above, the magnetic fields modulates the eigenfunction (something not seen by the absolute value) but otherwise the profile is remarkably stable. Theorem 1 tells us that the Filoche-Mayboroda landscape function (blue) remains an accurate predictor even as the magnetic field becomes active. Theorem 2 is illustrated by the red curve: the range of times (`$t$'s) used to compute the bound in Theorem 2 was $[0,10^{-4}]$ with time steps of size $10^{-6}$ and we see that the red curve is indeed an upper bound on the eigenfunction profile.
 
  \vspace{-10pt}
 \begin{center}
\begin{figure}[h!]
\includegraphics[width=0.7\textwidth]{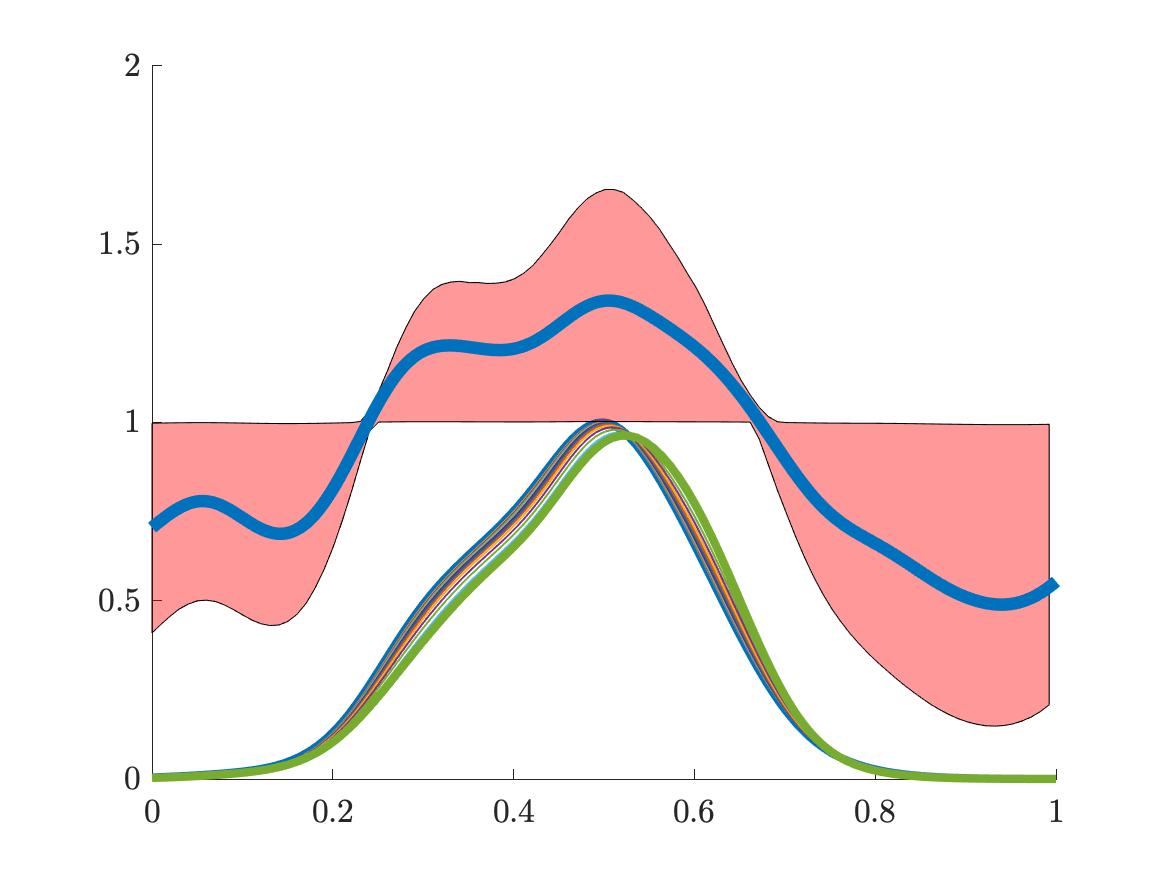}
\caption{Red shows the bounds coming from Theorem 2 for different values of $t$ which is compared to the Filoche-Mayboroda landscape function (in red). There are always values of $t$ for which it is better (the Corollary) and values $t$ for which it is worse.}
\end{figure}
\end{center}
 \vspace{-10pt}
 
We conclude with an illustation showing the role of the parameter $t$ in Theorem 2.  Fig. 5 shows, as the red region, a superposition of all the bounds that can be obtained for different values of $t$. As is easily seen from the closed-form expression of Theorem 2, very small values of $t$ simply recover the trivial bound $|\phi(x)|/\| \phi\|_{L^{\infty}} \leq 1$ corresponding to the straight line in Fig. 5. As $t$ grows, different types of bounds become available. What is particularly interesting is the contrast to the Filoche-Mayboroda landscape function (shown in blue). The main idea that underlies the proof of the Corollary is that the Filoche-Mayboroda landscape function arises as a weighted average of the bounds obtained by Theorem 2. Since it is a weighted average, there will be values $t$ for which the bound is at least as good (this is the content of the Corollary) and there will be other values where the bound is at most as good. This is nicely shown in Fig. 5: the Filoche-Mayboroda landscape function is sandwiched by the bounds obtained by Theorem 2.


\section{Proofs}
We first prove Theorem 2.  Theorem 1 follows immediately from this by averaging over time which then readily implies the Corollary.  Theorem 3 is proved via a different (though philosophically not unrelated) argument.
\subsection{Proof of Theorem 2}
\begin{proof} We abbreviate
$$ H(A,V) = \frac{1}{2} \left(- i\nabla - A(x)\right)^2 + V(x).$$
$H(A,V)$ is a self-adjoint operator with a real spectrum.  Observe that eigenfunctions
 $H(A,V) \psi = \lambda \psi$
have trivial behavior under the induced parabolic equation
$$ e^{-t H(A,V)} \psi = e^{-\lambda t} \psi.$$
The second part of the argument is based on rewriting the solution of the induced parabolic equation using
 the Feynman-Kac formula (see \cite{broderix})
$$ \left[e^{-t H(A,V)} \psi\right](x) = \mathbb{E}_{\omega}\left( e^{-S_t(A,V|\omega)} \chi_{\Omega}(\omega,t) \psi(\omega_x(t)) \right),$$
where the expectation is taken over all Brownian motion started at $x$, and where 
$$\chi_{\Omega}(\omega,t)  = \begin{cases} 1 \qquad &\mbox{if}~\forall ~0\leq s \leq t: \omega(s) \in \Omega \\ 0 \qquad &\mbox{otherwise} \end{cases}$$
measures whether the Brownian motion has left the domain. The quantity $S_t(A,V|\omega)$ is defined via the formula
\begin{align*}
 S_t(A,V|\omega) &= i \int_0^t A(\omega(s)) \cdot d\omega(s) + \frac{i}{2} \int_0^t (\nabla \cdot A)(\omega(s)) ds \\
 &+ \int_0^t V(\omega(s)) ds,
\end{align*}
where the first integral is an It\^{o} integral. An important aspect is that the first two terms are purely imaginary while the last term is purely real. This indicates
that the first two terms are responsible for oscillation while the last one is responsible for decay. After substituting the Feynman-Kac formula into the left-hand side equation  $e^{-t H(A,V)} \psi = e^{-\lambda t} \psi$ and multiplying by $e^{\lambda t},$ we see that
$$ \psi(x) = e^{\lambda t}  \left[e^{-t H(A,V)} \psi\right](x) = e^{\lambda t} \cdot  \mathbb{E}_{\omega}\left( e^{-S_t(A,V|\omega)} \chi_{\Omega}(\omega,t)\psi(\omega_x(t)) \right).$$
Let us now assume that $x_0 \in \Omega$ is arbitrary.
Then, using the fact that $S_t$ decomposes into real and imaginary parts,
the triangle inequality implies
\begin{align*}
  |\psi(x_0)| &\leq  e^{\lambda t} \cdot  \mathbb{E}_{\omega}\left(\chi_{\Omega}(\omega_{x_0}, t) \exp\left( - \int_0^t V(\omega_{x_0}(s)) ds \right) |\psi(\omega_{x_0}(t))| \right) \\
&\leq  \| \psi\|_{L^{\infty}}  \cdot e^{\lambda t} \cdot  \mathbb{E}_{\omega}\left(\chi_{\Omega}(\omega_{x_0}, t) \exp\left( - \int_0^t V(\omega_{x_0}(s)) ds \right) \right). 
 \end{align*}
Thus
$$ e^{-\lambda t} \cdot \frac{|\psi(x_0)|}{\|\psi\|_{L^{\infty}}} \leq   \mathbb{E}_{\omega}\left( \chi_{\Omega}(\omega_{x_0}, t) \exp\left( - \int_0^t V(\omega_{x_0}(s)) ds \right)  \right).$$
 \end{proof}
 
 \subsection{Proof of Theorem 1}
 \begin{proof}
Theorem 1 follows rather quickly from Theorem 2 by averaging in time. Recall that
 $$ e^{-\lambda t} \cdot \frac{|\psi(x_0)|}{\|\psi\|_{L^{\infty}}} \leq   \mathbb{E}_{\omega}\left( \chi_{\Omega}(\omega_{x_0}, t) \exp\left( - \int_0^t V(\omega_{x_0}(s)) ds \right)  \right).$$
 We will integrate both sides  over $0 \leq t \leq \infty$: the left-hand side is easy since
 $$ \int_0^{\infty} e^{-\lambda t} \cdot \frac{|\psi(x_0)|}{\|\psi\|_{L^{\infty}}}  = \frac{1}{\lambda}  \cdot \frac{|\psi(x_0)|}{\|\psi\|_{L^{\infty}}}.$$
 The right-hand side
\begin{align*}
X = \int_0^{\infty}  \mathbb{E}_{\omega}\left( \chi_{\Omega}(\omega_{x_0}, t) \exp\left( - \int_0^t V(\omega_{x_0}(s)) ds \right)  \right) dt
\end{align*}
requires some additional arguments. Note that by exchanging integration and expectation, we keep track of each Brownian particle until the exit time $\tau$, the  time where it exits the domain for the first time. Thus
\begin{align*}
X = \mathbb{E}_{\omega}  \int_0^{\tau}  \exp\left( - \int_0^{t} V(\omega_{x_0}(s)) ds \right)   d \omega.
\end{align*}
Finally, a quick inspection shows that this is merely the Feynman-Kac formula for $-\Delta +V$ with right-hand side 1 and therefore
$
X  = u(x_0).
$

 \end{proof}

 \subsection{Proof of the Corollary}
 \begin{proof} The Corollary is a simple consequence of the fact that we obtain Theorem 1 from Theorem 2 by averaging over time. We use the easy principle that any quantity has to be at least as large as its average (and also at most as large as it average) somewhere. More precisely, suppose now that the Corollary is incorrect: then, for all times $t > 0$, we have
$$  \mathbb{E}_{\omega}\left( \chi_{\Omega}(\omega_{x_0}, t) \exp\left(\lambda \cdot t - \int_0^t V(\omega_{x_0}(s)) ds \right)  \right) > \lambda \cdot u(x)$$
which can be rewritten as
$$  \mathbb{E}_{\omega}\left( \chi_{\Omega}(\omega_{x_0}, t) \exp\left( - \int_0^t V(\omega_{x_0}(s)) ds \right)  \right) > \lambda e^{-\lambda t} \cdot u(x).$$
Integrating again both sides for $0 \leq t \leq \infty$ in the same way as in the proof of Theorem 1 now leads to a contradiction.
\end{proof}

 \subsection{Proof of Theorem 3}
 \begin{proof}
 We start with the probabilistic representation of the landscape function already used in the proof of Theorem 1
 \begin{align*}
  u(x) &= \mathbb{E}_{\omega} \int_0^{\tau} \exp \left( - \int_0^s V(\omega_x(t)) dt\right) ds.
  \end{align*}
 where $\tau$ is the exit time of the domain. 
 We will use this formula with a small variation: the Dirichlet condition on the solution requires us to stop the integral at the exit time $\tau$. 
 For convenience, we will assume in this proof that the particle does not stop contributing but merely gets stuck in place and keeps contributing allowing us to write, tautologically, 
  \begin{align*}
u(x) &= \mathbb{E} \int_0^{\tau} \exp \left( - \int_0^s V(\omega_x(t)) dt\right) ds\\
&=   \mathbb{E} \int_0^{\infty}   \exp \left( - \int_0^{s } V(\omega_x(t)) dt\right) ds -  \mathbb{E} \int_{\tau}^{\infty}   \exp \left( - \int_0^{s } V(\omega_x(t)) dt\right) ds.
  \end{align*}
  Note all integrals are finite since $V \geq V_{\max} > 0$.
We shall now bound the first term from below and the second term from above.

  For the first term, Jensen's inequality implies
  \begin{align*}
 \mathbb{E} \int_0^{\infty} \exp \left( - \int_0^s V(\omega_x(t)) dt\right) ds   &= \int_0^{\infty}   \mathbb{E} \exp \left( - \int_0^s V(\omega_x(t)) dt\right) ds \\
 &\geq \int_0^{\infty} \exp \left( -  \mathbb{E}  \int_0^s V(\omega_x(t)) dt\right) ds \\
   &=   \int_0^{\infty} \exp \left( - s \cdot  \mathbb{E} \frac{1}{s} \int_0^s V(\omega_x(t)) dt\right) ds \\
&\geq     \int_0^{\infty} \exp \left( - s \sup_{w > 0}  \mathbb{E} \frac{1}{w} \int_0^w V(\omega_x(t)) dt\right) ds \\
&= 1/\left( \sup_{w > 0}  \mathbb{E} \frac{1}{w} \int_0^w V(\omega_x(t)) dt \right).
   \end{align*}
   It remains to analyze the second term.
 We note that $V \geq V_{\min} > 0$ and thus
 \begin{align*}
   \mathbb{E} \int_{\tau}^{\infty}   \exp \left( - \int_0^{s } V(\omega_x(t)) dt\right) ds &\leq   \mathbb{E} \int_{\tau}^{\infty}   \exp \left( - s V_{\min}\right) ds \\
   &= \frac{\mathbb{E} e^{-\tau V_{\min}}}{V_{\min}},
   \end{align*}
   where, as mentioned above, $\tau$ is the exit time for a Brownian motion started at $x.$
   It is therefore sufficient to find a lower bound on $\mathbb{E} e^{-\tau V_{\min}}$.
   
 Here, we use a basic monotonicity property. If $x \in \Omega \subseteq \Omega_2$, then $\tau_{\Omega} \leq \tau_{\Omega_2}$ since any particle exiting $\Omega_2$ must also have exited $\Omega$. Thus
 $$ \mathbb{E} ~e^{-\tau_{\Omega_2} V_{\min}} \leq  \mathbb{E}~ e^{-\tau_{\Omega} V_{\min}}.$$
  We use this by considering the largest ball centered at $x$ and fully contained inside the domain, $B_r(x) \subseteq \Omega \subset \mathbb{R}^d$. Tautologically, $r = d(x, \partial \Omega)$. This in turn is a relatively well-studied quantity \cite{hamana}.
  Setting 
  $$ \nu = \frac{d}{2}-1,$$
  we have for all $\lambda > 0$ that 
  $$ \mathbb{E} \left( e^{- \lambda \cdot \tau_{B_r(x)}}\right) = \frac{( r \sqrt{2 \lambda})^{\nu}}{2^{\nu} \Gamma(\nu + 1)} \frac{1}{I_{\nu}(r \sqrt{2\lambda})}.$$
  Given the explicit form of the expression and $\lambda \geq 0$, we clearly have
  $$  \mathbb{E} \left( e^{- \lambda \cdot \tau_{B_r(x)}}\right) \leq 1$$
  and the expression converges to 1 as $r \rightarrow 0$. For large values of $r$, we may use the standard asymptotic result
  $$ I_{\nu}(x) = (1+o(1)) \frac{e^x}{\sqrt{2\pi x}}$$
  to deduce that
    $$ \mathbb{E} \left( e^{- V_{\min} \cdot \tau_{B_r(x)}}\right) \leq c_{\nu} \cdot \exp\left(- r \cdot \sqrt{V_{\min}}\right).$$
   
 \end{proof}

\end{document}